\documentclass[12pt]{article}

\usepackage{amsmath}
\usepackage{amssymb}
\usepackage{amsthm}
\usepackage{a4wide}
\usepackage[english]{babel}

\numberwithin{equation}{section}

\newtheorem{Th}{Theorem}
\newtheorem{Lemma}{Lemma}
\newtheorem{Prop}{Proposition}
\newtheorem{Cor}{Corollary}

\newcommand{\Dem}{\noindent{\bf Proof }}

\newcommand{\N}{\mathbb{N}}
\newcommand{\Z}{\mathbb{Z}}
\newcommand{\Q}{\mathbb{Q}}
\newcommand{\R}{\mathbb{R}}

\newcommand{\eps}{\varepsilon}

\newcommand{\lcm}{{\rm lcm}}

\newcommand{\nve}{{\varepsilon_1}}

\newcommand{\ffi}{\varphi}
\newcommand{\nvffi}{\ffi}

\newcommand{\Span}{{\rm Span}}

\newcommand{\lam}{\lambda}

\newcommand{\capasig}{\kappa_\sigma }
\newcommand{\capasigpr}{\kappa_{\sigma'}}
\newcommand{\capaId}{\kappa_{\Id}}
\newcommand{\etasig}{\eta_\sigma}
\newcommand{\etasigpr}{\eta_{\sigma'}}
\newcommand{\etaId}{\eta_{\Id}}

\newcommand{\unk}{\{1,\ldots,k\}}
\newcommand{\capasigpuk}{\{\capasig+1,\ldots,k\}}
\newcommand{\capasigk}{\{\capasig,\ldots,k\}}
\newcommand{\unp}{\{1,\ldots,p\}}
\newcommand{\unpmu}{\{1,\ldots,p-1\}}

\newcommand{\moins}{\setminus}

\newcommand{\Netoile}{\mathbb{N}^*}

\newcommand{\sk}{{\mathfrak S}_k}
\newcommand{\Id}{{\rm Id}}

\newcommand{\nvcalC}{{\mathcal{C}(\eps,Q)}}
\newcommand{\nvLam}{\Lambda(Q)}
\newcommand{\nvLamlz}{\Lambda(Q_{\ell_0})}

\newcommand{\calC}{{\mathcal{C}}}

\newcommand{\combitiny}[2]{{\tiny \left( \! \!  \begin{array}{c} #1 \\ #2 \end{array} \! \! \right)}}
\newcommand{\combi}[2]{{  \left( \begin{array}{c} #1 \\ #2 \end{array} \right)}}

\newcommand{\pha}{\varphi}

\renewcommand{\Im}{{\rm Im}}

\newcommand{\rk}{{\rm rk}}

\newcommand{\Dist}{{\rm Dist}}

\title{Nesterenko's linear independence criterion for vectors}
\author{St\'ephane Fischler}
\date{\today}

\begin{document}

\newcommand{\pb}[1]{{\bf #1}}

\newcommand{\Li}{{\rm Li}}
\newcommand{\norm}[1]{\lVert #1\rVert}

\newcommand{\om}{\omega}
\newcommand{\unkmu}{\{1,\ldots,k-1\}}

\maketitle 

\begin{abstract}
In this paper we deduce a lower bound for the rank of a family of $p$ vectors in $\R^k$ (considered as a vector space over the rationals) from the existence of a sequence of linear forms on $\R^p$, with integer coefficients, which are small at $k$ points. This is a generalization to vectors of Nesterenko's linear independence criterion (which corresponds to $k=1$), used by Ball-Rivoal to prove that infinitely many values of Riemann zeta function at odd integers are irrational. The proof is based on geometry of numbers, namely Minkowski's theorem on convex bodies.
\end{abstract}

\section{Introduction}

The motivation for this paper comes from irrationality results on values of Riemann zeta function $\zeta(s) = \sum_{n=1}^\infty \frac1{n^s} $ at odd integers $s \geq 3$. The first result is due to Ap\'ery \cite{Apery}: $\zeta(3) \not\in \Q$. The next breakthrough in this topic is due to Rivoal \cite{RivoalCRAS} and Ball-Rivoal \cite{BR}: 
\begin{equation} \label{eqBR}
\dim_\Q \Span_\Q (1 , \zeta(3),  \zeta(5),  \zeta(7), \ldots,  \zeta(a)) \geq \frac{\log a}{1+\log 2 } (1+o(1))
\end{equation}
as $ a\to\infty$, where $a$ is an odd integer; notice this is a lower bound on the rank of this family of real numbers, in $\R$ considered as a vector space over the rationals. Conjecturally the left handside is equal to $\frac{a+1}2$, but even the constant $\frac1{1+\log 2}$ in Eq. \eqref{eqBR} has never been improved. Actually,  known refinements of Ball-Rivoal's proof provide sharper lower bounds only for fixed values of $a$: the improvement always lies inside the error term $o(1)$ as $a \to \infty$. 

\medskip

However, the following improvement of \eqref{eqBR} is proved in \cite{SFdistrib}:

\begin{Th} \label{thdistrib}
Let $\eps >0$, and $a$ be an odd integer sufficiently large with respect to $\eps$. Then letting $N$ denote the integer part of $\frac{1-\eps}{1+\log 2}\log a$, there exist odd integers $\sigma_1,\ldots,\sigma_N$ between 3 and $a$ such that:
\begin{itemize}
\item 1, $\zeta(\sigma_1)$, \ldots, $\zeta(\sigma_N)$ are linearly independent over the rationals;
\item For any $i\neq j$, $|\sigma_i  - \sigma_j| > a^\eps$.
\end{itemize}
\end{Th}

In particular, if there are only $N$ odd integers $\sigma$ between 3 and $a$ such that $\zeta(\sigma)$ is irrational, then they have to be evenly distributed (see  \cite{SFdistrib}).

\bigskip

The strategy for proving Theorem \ref{thdistrib} is based on the following classical construction. For non-negative integers $\beta$, $b$, $n$, $r$ with $\beta $ and $b$ odd, $1\leq \beta \leq b$, and $2br  < a$, let 
\begin{equation} \label{eqdefibnnestsev}
J_{\beta,n} =  \frac{ d_{2n}^{a+b-1} (2n)!^{a-2b r } }{( \beta-1)! }  \sum_{k =  1}^\infty  \frac{{\rm d}^{\beta-1}}{{\rm d}k ^{\beta-1}} \Big( \frac{(k- 2r n)_{2rn}^b (k+2n+1)_{2rn}^b }{(k)_{2n+1}^a} \Big),
\end{equation}
where the derivative is taken at $k$, Pochhammer's symbol is defined by $(\alpha)_p = \alpha (\alpha+1)\ldots(\alpha+p-1)$, and $d_{2n}$ is the least common multiple of 1, 2, 3, \ldots, $2n$. It is not difficult   to prove that 
$$J_{\beta,n}  = \widetilde{\ell}_{\beta,n}  +   \ell_{3,n}   \combi{\beta+1}{\beta-1} \zeta(\beta+2) +  \ell_{5,n}  \combi{\beta+3}{\beta-1}  \zeta(\beta+4) + \ldots + \ell_{a,n}   \combi{\beta+a-2}{\beta-1} \zeta(\beta+a-1)$$
with  integers $ \widetilde{\ell}_{\beta,n}$ and $\ell_{i,n}$; moreover $J_{\beta,n}  $ tends to 0 as $n\to\infty$, for any $\beta$, provided the parameters satisfy suitable relations (and up to technicalities, see  \cite{SFdistrib} for precise statements). This   can be seen as a sequence $(L_n)$ of linear forms on $\R^{(a+b)/2}$,  with integer coefficients, that take small values $L_n(e_j) = J_{2j-1,n}$ at $ k = \frac{b+1}2$ points $e_1,\ldots,e_k \in \R^{(a+b)/2}$. The key point in the proof of Theorem \ref{thdistrib} is then to apply the following result, to which the present paper is devoted.

\bigskip

 We let  $\R^p$ be endowed with its canonical scalar product  and the corresponding norm.

\bigskip

\begin{Th} \label{thintro2}
Let $1\leq k \leq p-1$,  and   $e_1,\ldots,e_k \in \R^p$. 

Let $\tau_1,\ldots,\tau_k >0 $ be pairwise distinct  real numbers. 

Let $(Q_n)_{n \geq 1}$  be an increasing sequence of positive integers, such that   $Q_{n+1} = Q_n^{1+o(1)}$.

For any $n\geq 1$, let $ L_n = \ell_{1,n}X_1 + \ldots + \ell_{p,n}X_p$  be a linear form on $\R^p$, with integer coefficients $\ell_{i,n}$ such  that, as $n\to\infty$:
$$|L_n(e_j) | = Q_n^{-\tau_j+o(1)}    \mbox{ for any } j \in \unk \mbox{ and }  
 \max_{1 \leq i \leq p} |\ell_{i,n} | \leq Q_n^{1+o(1)}.
$$

Then:
\begin{enumerate}
\item[$(i)$] If  $F$ is a subspace of $\R^p$  defined over $\Q$ which contains $e_1$, \ldots, $e_k$ then 
$$  \dim F \geq k + \tau_1 +\ldots + \tau_k .$$
In other words, letting $C_1 , \ldots,  C_p \in \R^k$ denote the columns of the matrix whose rows are $e_1,\ldots,e_k \in\R^p$, we have
$$\rk_\Q(C_1,\ldots,C_p)\geq k + \tau_1 +\ldots + \tau_k  $$
in $\R^k$ seen as a $\Q$-vector space.
\item[$(ii)$] The vectors  $e_1,\ldots,e_k $ are $\R$-linearly independent in $\R^p$, and the $\R$-subspace they span  does not intersect $\Q^p \moins \{(0,\ldots,0)\}$.
\item[$(iii)$] Let $\eps>0$, and $Q$ be sufficiently large (in terms of $\eps$). Let $\nvcalC$ denote the set of all vectors that can be written as $\lambda_1 e_1 + \ldots + \lambda_k e_k  + u$ with:
$$\left\{ \begin{array}{l}
 \lambda_1,\ldots,\lambda_k\in \R  \mbox{  such that } |\lambda_j| \leq Q^{\tau_j-\eps}   \mbox{ for any } j\in\unk\\
 u \in (\Span_\R(e_1,\ldots,e_k))^\perp \mbox{ such that } \norm{u} \leq Q^{-1-\eps} 
\end{array}\right.
$$
Then $\nvcalC \cap \Z^p = \{(0,\ldots,0)\}$.
\end{enumerate}
\end{Th}

If $k=1 $  this is exactly Nesterenko's linear independence criterion \cite{Nesterenkocritere} used in the proof of Ball-Rivoal's result \eqref{eqBR}.

\bigskip

In the conclusions, $(ii)$ is an easy result, and $(iii)$ is the main part (it is a quantitative version of $(ii)$). We deduce $(i)$ from $(iii)$ using Minkowski's convex body theorem, thereby generalizing the proof given in \cite{SFZu} and \cite{eddzero} of  Nesterenko's linear independence criterion. The equivalence between both statements of $(i)$ comes from linear algebra; it is proved in \S \ref{subsecrank}.

\bigskip

A result analogous to  Theorem \ref{thintro2}, but in which $p$ linearly independent  linear forms like $L_n$ appear in the assumption, is proved in \S \ref{subsecsiegel}. This  linear independence criterion (in the style of Siegel's) is much easier to prove than Theorem \ref{thintro2}. Both results can be thought of as transference principles. In this respect it is worth pointing out that in Theorem \ref{thintro2} we assume essentially that for any positive integer $Q$ there is a linear form : indeed this is $L_n$, where $n$ is such that   $Q_n \leq Q < Q_{n+1}$ so that $Q = Q_n^{1+o(1)}$ because    $Q_{n+1} = Q_n^{1+o(1)}$. The assumptions imply that this linear form belongs to some convex body, and conclusion $(iii)$ asserts that (up to $Q^\eps$) the dual convex body does not contain any non-zero integer point. Therefore it is reasonable to imagine that $(iii)$ is an optimal conclusion up to $Q^\eps$. In general  the lower bound $k + \tau_1 +\ldots + \tau_k $ in $(i)$ is optimal too (see \cite{FHKL} for a converse statement, valid almost everywhere). In the special case  $p=2$, $k=1$, and $e_1=(1,\xi)$, Theorem  \ref{thintro2}  $(iii)$ yields an upper bound $\mu(\xi) \leq 1 + \frac1{\tau_1}$ on the irrationality exponent of $\xi$, and reduces essentially to Lemma 1  of \cite{eddzero}. A converse statement in this case is proved in \cite{eddzero} (Theorem 1).

\bigskip

The assumption that $\tau_1$, \ldots, $\tau_k$ are pairwise distinct is very important in Theorem \ref{thintro2}, and it cannot be omitted. For instance, if $\tau_1=\tau_2$ then  $L_n(e_1-e_2)$ could be very small: up to replacing $(e_1,e_2)$ with $(e_1+e_2,e_1-e_2)$, this amounts to dropping the assumption that the linear forms $L_n$ are not too small at the points $e_j$. Now this assumption is known to be essential, already in the classical case of Nesterenko's linear independence criterion (except for proving the linear independence of   three numbers, see Theorem 2  of \cite{SFZu}). Actually, if  $\tau_1=\tau_2$ then  $L_n(e_1-e_2)$ could even vanish, so the possibility that $e_1=e_2$ cannot be eliminated: even assertion $(ii)$ may fail to hold.

\bigskip

We shall prove Theorem \ref{thintro2} in a more general form, stated in \S \ref{secenonce}, which allows the sequences $(|L_n(e_j)|)_{n\geq 1}$ to oscillate (as in \cite{SFoscillate}), and  takes into account divisors of the coefficients $\ell_{i,n}$ (as in \cite{SFZu}); the former is used in \cite{SFdistrib} to prove Theorem \ref{thdistrib}. We also include a refinement useful when $L_n$ is not too large at some other point, which is new even in the classical case of Nesterenko's linear independence (with $k=1$).

\bigskip

We hope that our results will have Diophantine applications besides those of \cite{SFdistrib}; we mention some directions in \S \ref{secdio}, connected to polylogarithms or zeta values. Our criterion could be used also for $q$-analogues, as in \cite{SFZu}.

\bigskip

The structure of this text is as follows. In \S \ref{secenonce} we state our result in a very general form, of which Theorem \ref{thintro2} is a special case. Section \ref{sec3} is devoted to the proof; then we deduce some corollaries in \S \S \ref{subsec21} and \ref{subsec22}. We prove an analogous result in the style of Siegel's linear independence criterion in \S \ref{subsecsiegel}, and conclude in \S \ref{secdio} with Diophantine applications.

\section{Statement of the criterion} \label{secenonce}

The following generalization of Theorem \ref{thintro2} is our main result.

\begin{Th} \label{thgal}
Let $1\leq k \leq p-1$, and $e_1,\ldots,e_k\in\R^p$. Let $(v_1,\ldots,v_p)$ denote a basis of $\R^p$. Let $\tau_1,\ldots,\tau_k>0$, $\sigma_1\geq \ldots \geq \sigma_p>0$,  $\om_1,\ldots,\om_k$, $\pha_1,\ldots, \pha_k$ be real numbers, with $\tau_1,\ldots,\tau_k$ pairwise distinct. Assume that there exist infinitely many integers $n$ with the following property: for any $j \in \unk$, $n \om_j + \pha_j \not\equiv \frac{\pi}2 \bmod \pi$.

Let $(Q_n)_{n \geq 1}$  be an increasing sequence of positive integers, such that  $Q_{n+1} = Q_n^{1+O(1/n)}$; if $\om_1=\ldots=\om_k=0$, this assumption can be weakened to
$Q_{n+1} = Q_n^{1+o(1)}$.

For any $n\geq 1$, let $ L_n = \ell_{1,n}X_1 + \ldots + \ell_{p,n}X_p$  be a linear form on $\R^p$, with integer coefficients $\ell_{i,n}$ such  that, as $n\to\infty$:
\begin{equation} \label{eqhypasyth}
|L_n(e_j) | = Q_n^{-\tau_j+o(1)} |\cos(n \om_j + \pha_j) + o(1)|  \mbox{ for any } j \in \unk,
\end{equation}
and 
$$|L_n(v_i)|\leq Q_n^{\sigma_i+o(1)}\mbox{ for  any } i\in\unp.$$
For all $n\geq 1 $ and $i\in\unp$, let $\delta_{i,n}$ be a positive divisor of $\ell_{i,n}$ such that:
\begin{itemize}
\item[$(i)$] $\delta_{i,n}$ divides $\delta_{i+1,n}$ for any $n\geq 1$ and any $i\in\unpmu$,
\item[$(ii)$] $\frac{\delta_{j,n}}{\delta_{i,n}}$ divides $\frac{\delta_{j,n+1}}{\delta_{i,n+1}}$ for any $n\geq 1$ and any $0 \leq i < j \leq p$, with $\delta_{0,n} = 1$,
 \item[$(iii)$] $\delta_{i,n} = Q_n^{d_i+o(1)}$ as $n\to\infty$ for any $i\in\unp$, with real numbers $d_i$ such that $0 \leq d_1 \leq \ldots\leq d_p\leq \sigma_p$.
\end{itemize}

Then:
\begin{enumerate}
\item[$(i)$] If  $F$ is a subspace of $\R^p$  defined over $\Q$ which contains $e_1$, \ldots, $e_k$ then $s = \dim F$ satisfies $s \geq k+1$ and 
\begin{equation} \label{eqqltthgal}
 \sigma_1 + \ldots + \sigma_{s-k} \geq   \tau_1 +\ldots + \tau_k  + d_1 + \ldots + d_s.
 \end{equation} 
In other words, letting $C_1 , \ldots,  C_p \in \R^k$ denote the columns of the matrix whose rows are $e_1,\ldots,e_k \in\R^p$,  the rank $s$ of the family $(C_1,\ldots,C_p)$ in $\R^k$  seen as a $\Q$-vector space satisfies $s \geq k+1$ and  Eq. \eqref{eqqltthgal}.
\item[$(ii)$] The vectors  $e_1,\ldots,e_k $ are $\R$-linearly independent in $\R^p$, and the $\R$-subspace they span  does not intersect $\Q^p \moins \{(0,\ldots,0)\}$.
\item[$(iii)$] Let $\eps>0$, and $Q$ be sufficiently large (in terms of $\eps$). Let $\nvcalC$ denote the set of all vectors that can be written as $\lambda_1 e_1 + \ldots + \lambda_k e_k  + u$ with:
$$\left\{ \begin{array}{l}
 \lambda_1,\ldots,\lambda_k\in \R  \mbox{  such that } |\lambda_j| \leq Q^{\tau_j-\eps}   \mbox{ for any } j\in\unk\\
 u \in (\Span_\R(e_1,\ldots,e_k))^\perp \mbox{ such that }  u = \mu_1v_1+\ldots+\mu_pv_p \mbox{ with } |\mu_i | \leq Q^{-\sigma_i-\eps}\\
 \hspace{10cm} \mbox{ for any } i\in\unp.
\end{array}\right.
$$
Let $\nvLam$ denote the set of all $(x_1,\ldots,x_p) \in\Q^p$ such that $\delta_{i,\Psi(Q)}x_i\in\Z$ for any $i\in\unp$, where $\Psi(Q)$ is the largest integer $n$ such that $Q_n\leq Q$. 

Then $\nvcalC \cap \nvLam = \{(0,\ldots,0)\}$.
\end{enumerate}
\end{Th}

In the special case where $\sigma_i = \delta_{i,n}=1$ and $d_i = \om_j = \pha_j=0$ for any $i$, $j$, $n$, and $(v_1,\ldots,v_p)$ is the canonical basis of $\R^p$, this is exactly Theorem \ref{thintro2} stated in the introduction. Indeed Eq. \eqref{eqhypasyth} reads $|L_n(e_j)| = Q_n^{-\tau_j+o(1)}$ in this case, and we have $L_n(v_i) = \ell_{i,n}$; moreover Eq. \eqref{eqqltthgal} reads 
$$\dim F \geq k + \tau_1+\ldots+\tau_k.$$
There is only a minor difference in $(iii)$, where the norm of $u$ is the Euclidean one in Theorem \ref{thintro2}, and the infinite one  in Theorem \ref{thgal}; of course this is not significant.

\bigskip

The real numbers $\om_j$ and $\pha_j$ allow oscillating behaviors of the sequences $(|L_n(e_j)|)_{n\geq 1}$. This is used in \cite{SFdistrib}, where the saddle point method is applied. In the special case of Theorem \ref{thintro2} with $k=1$, the corresponding generalization of Nesterenko's linear independence criterion has been proved in \cite{SFoscillate} when $Q_n= \beta^n$ for some $\beta>1$ (which is the most interesting case). We generalize it here to any sequence $(Q_n)$ such that $Q_{n+1}  = Q_n^{1+O(1/n)}$; eventhough this assumption is slightly more restrictive than  the usual one $Q_{n+1}  = Q_n^{1+o(1)}$, it is still general enough to include sequences $Q_n = \beta^{n^d}$ with $\beta>1$ and $d>0$. 

\bigskip

The divisors $\delta_{i,n}$ allow one to make use of divisibility properties of the coefficients $\ell_{i,n}$: for instance, in most  constructions of linear forms in zeta values, $\ell_{i,n}$ is a multiple of $\delta_{i,n} = d_n^{e_i}$ for some $e_i\geq 1$, where $d_n = \lcm(1,2,\ldots,n)$.  These divisors are used in \cite{SFdistrib} to prove a variant of Theorem \ref{thdistrib}. The first refinement of Nesterenko's linear independence criterion involving such divisors $\delta_{i,n}$ is Theorem 1 of \cite{SFZu}, which is essentially the special case of Theorem \ref{thgal} $(i)$ where $k=1$, $\sigma_i=1$, $\om_j=\pha_j=0$, and  $(v_1,\ldots,v_p)$ is the canonical basis of $\R^p$; it is the main ingredient in the proof \cite{SFZu} that 1, $\zeta(3)$ and $\zeta(j)$ are $\Q$-linearly independent for some odd integer $j$ between 5 and 139.

\bigskip

The real numbers $\sigma_i$ allow one to take advantage of the fact that the linear forms $L_n$ might be smaller than $\norm{L_n}$ at some given points $v_i$ (eventhough $L_n(v_i)$ does not tend to 0 as $n\to\infty$). For instance, if  $(v_1,\ldots,v_p)$ is the canonical basis, this is useful when one has a sharper upper bound on $|\ell_{i,n}|$ for some values of $i$ than for others. This feature is new even in the case of Nesterenko's linear independence criterion  (namely, with  $k=1$, $\sigma_i = \delta_{i,n}=1$, and $d_i = \om_j = \pha_j=0$). It would be interesting to deduce from this refinement a Diophantine consequence. Actually it happens   for linear forms in zeta values that  $\lim_{n\to\infty} |\ell_{i,n}|^{1/n}$ exists for any $i$ and does depend on $i$. For instance, F. Amoroso and T. Rivoal have noticed that 
in the expansion  of 
$$n!^{a-1} \sum_{k=1}^\infty \frac{(k-n)_n}{(k)_{n+1}^a}$$
as a linear combination of zeta values, the coefficients of odd and even zeta values don't have the same size (provided $a$ is even).

\bigskip

It is very important in Theorem \ref{thgal} that $\tau_1,\ldots,\tau_k$ are pairwise distinct; however it is not always necessary to compute their exact values. For instance, if $\min(\tau_1,\ldots,\tau_k)$ is greater than or equal to some $\tau>0$, then Eq. \eqref{eqqltthgal} implies
$$\sigma_1+\ldots+\sigma_{s-k}\geq k\tau + d_1+\ldots+d_s;$$
in the special case of Theorem \ref{thintro2} this lower bound reads $\dim F \geq k(1+\tau)$. This remark is already used (with $k=1$) in \cite{BR}, and also in the proof \cite{SFdistrib} of Theorem \ref{thdistrib}. We refer to \S \ref{subsec22} below for a related result.

\bigskip

At last,  notice that if the assumptions of Theorem \ref{thgal} hold with $e_1,\ldots,e_k$, then they hold also if we forget one of the $e_j$'s (say $e_k$, with $k\geq 2$). The same implication holds also for parts $(ii)$ and $(iii)$ of the conclusion, since the convex body $\nvcalC$ becomes smaller when $e_k$ is omitted. However this implication does not hold for part $(i)$; to fix this we refine part $(i)$ as follows.

\begin{Cor} \label{corgal}
In the situation of Theorem \ref{thgal}, assume also that $\tau_1 > \ldots > \tau_k$. Then for any  subspace $F$ of $\R^p$  defined over $\Q$  we have
\begin{equation} \label{eqcorgal}
s \geq t+1 \mbox{ and }  \sigma_1 + \ldots + \sigma_{s-t} \geq    \tau_{k+1-t}  +\ldots + \tau_k  + d_1 + \ldots + d_s, 
 \end{equation} 
provided that $s = \dim F$ and $t = \dim( F\cap \Span_\R(e_1,\ldots,e_k))$ are positive.

 In other words, for any surjective $\R$-linear map $\pi : \R^k \to \R^t $ with $t\geq 1$, Eq. \eqref{eqcorgal} holds with
 $$s = \rk_\Q(\pi(C_1),\ldots,\pi(C_p))$$
 where the rank is computed in $\R^t$ seen as a $\Q$-vector space. 
\end{Cor}

\Dem of Corollary \ref{corgal}: Let   $F$ be a subspace of $\R^p$  defined over $\Q$; assume that $s = \dim F$ and $t = \dim( F\cap \Span_\R(e_1,\ldots,e_k))$ are positive. For any $j\in\unk$ we let $D_j = \dim( F\cap \Span_\R(e_1,\ldots,e_j))$, so that $0 \leq D_1\leq\ldots\leq D_k=t$ and $D_j \in\{D_{j-1}, D_{j-1}+1\}$ for any $j$ (with $D_0=0$). Then there exist $t$ integers $1\leq j_1 < \ldots < j_t\leq k$ such that $D_j = D_{j-1}+1$ if, and only if, $j$ is among the $j_i$'s. For any $i\in\{1,\ldots,t\}$, there exists  $e'_i\in F\cap \Span_\R(e_1,\ldots,e_{j_i}) $ such that $e'_i\not\in   \Span_\R(e_1,\ldots,e_{j_i-1}) $. Then we have $e'_i = \sum_{j=1}^{j_i} \lambda_{i,j} e_j$ for real numbers $\lambda_{i,j}$ such that $\lambda_{i, j_i}\neq 0$. Since $\tau_1> \ldots > \tau_{j_i}$, Eq. \eqref{eqhypasyth} yields
$$| L_n(e'_i)| = Q_n^{-\tau_{j_i}+o(1)}  |\cos(n \om_{j_i} + \pha_{j_i}) + o(1)| .$$
Therefore Theorem \ref{thgal} applies to $e'_1,\ldots, e'_t$ with $\tau_{j_1}, \ldots, \tau_{j_t}$. Since $\tau_1> \ldots > \tau_k$, the inequality \eqref{eqqltthgal} obtained in this way implies Eq. \eqref{eqcorgal}. This concludes the proof of Corollary \ref{corgal}, except for the second part of the conclusion which will be proved at the end of \S \ref{subsecrank} below.

\section{Proof of the criterion} \label{sec3}

This section is devoted to proving  Theorem \ref{thgal}, of which   Theorem \ref{thintro2} stated in the introduction is a special case (see \S \ref{secenonce}).  Reindexing $e_1$, \ldots, $e_k$ is necessary, we assume that $\tau_1 > \ldots > \tau_k > 0$. This assumption  will be used in \S\S \ref{subsecdemii} and \ref{subsec33}.

\subsection{Rational rank of vectors} \label{subsecrank}
 
 In this section, we give some details about the conclusions of our criterion, which allow us to  prove the equivalence of both conclusions of $(i)$ in Theorems \ref{thintro2} and \ref{thgal}, and to conclude the proof of Corollary \ref{corgal}.

\bigskip

In Nesterenko's linear independence criterion, a lower bound is derived for the dimension of the $\Q$-subspace of $\R$ spanned by $\xi_0,\ldots,\xi_r\in \R$, that is, for the $\Q$-rank of $\xi_0,\ldots,\xi_r$ in $\R$ considered as a vector space over $\Q$. This rank is equal to the dimension of the smallest subspace of $\R^{r+1}$, defined over the rationals, which contains the point $(\xi_0,\ldots,\xi_r)$. We generalize in Lemma \ref{lemegal} below this equality to our setting.

Recall that a subspace $F$ of $\R^p$ is said to be {\em defined over $\Q$} if it is the zero locus of a family of linear forms with rational coefficients. This is equivalent to the existence of a basis (or a generating family) of $F$, as a vector space over $\R$, consisting in vectors of $\Q^p$ (see for instance \S 8 of \cite{Boualg}).  Since the intersection of a family of subspaces of $\R^p$ defined over $\Q$ is again defined over $\Q$, there exists for any subset $S \subset \R^p$ a minimal subspace of $\R^p$, defined over $\Q$, which contains $S$: this is the intersection of all  subspaces of $\R^p$, defined over $\Q$, which contain $S$. 

Let $M$ be a matrix with $k \geq 1$ rows, $p \geq 1$ columns, and real entries. Letting $e_1,\ldots,e_k \in \R^p$ denote the rows of $M$, we can consider as above the smallest subspace of $\R^p$, defined over $\Q$, which contains  $e_1,\ldots,e_k $. On the other hand, we denote by $C_1,\ldots,C_p \in \R^k$ the columns of $M$ and consider $\R^k$ as an infinite-dimensional vector space over $\Q$. Then $\Span_\Q (C_1,\ldots,C_p)$ is the smallest $\Q$-vector subspace of $\R^k$ containing $C_1,\ldots,C_p$; it consists in all linear combinations $r_1C_1+\ldots+r_p C_p$ with $r_1,\ldots,r_p\in \Q$. Its dimension (as a $\Q$-vector space) is the rank (over $\Q$) of $C_1,\ldots,C_p$, denoted by $\rk_\Q (C_1,\ldots,C_p)$. 

\begin{Lemma} \label{lemegal} 
Let $M \in {\rm Mat}_{k,p}(\R)$ with $k,p \geq 1$. Denote by $e_1,\ldots,e_k \in \R^p$ denote the rows of $M$, and by $C_1,\ldots,C_p \in \R^k$ its columns. Then $\rk_\Q (C_1,\ldots,C_p)$ is the dimension of the smallest subspace of $\R^p$, defined over $\Q$, which contains  $e_1,\ldots,e_k $. 
\end{Lemma} 

When $k=1$, this lemma means that the $\Q$-rank of $\xi_0,\ldots,\xi_r$   is equal to the dimension of the smallest subspace of $\R^{r+1}$, defined over the rationals, which contains the point $(\xi_0,\ldots,\xi_r)$. 

\bigskip

\Dem of Lemma \ref{lemegal}: Let $G = (\Span_\R(e_1,\ldots,e_k ))^\perp$, where $\R^p$ is equipped with the usual scalar product. Let $F$ denote the minimal  subspace of $\R^p$, defined over $\Q$, which contains  $e_1,\ldots,e_k $.  Then $F^\perp$ is the maximal subspace of $\R^p$, defined over $\Q$, which  is contained in $G =  \{e_1,\ldots,e_k \}^\perp$. Therefore $F^\perp = \Span_\R(G \cap \Q^p)  = (G\cap \Q^p)  \otimes_\Q \R$: any basis of the $\Q$-vector space $G \cap\Q^p$ is an $\R$-basis of $F^\perp$. Since $G \cap \Q^p = \ker \psi$ where $\psi : \Q^p \to \R^k$ is defined by $\psi(r_1,\ldots,r_p) = 
r_1C_1+\ldots+r_p C_p$, we have:
$$\dim_\R F  = p - \dim_\R F^\perp = p - \dim_\Q (G \cap\Q^p)  = \rk_\Q \psi = \rk_\Q (C_1,\ldots,C_p).$$

This concludes the proof of Lemma \ref{lemegal}.

\bigskip

Let us deduce from Lemma \ref{lemegal} the following generalization, and use it to prove the second assertion of Corollary \ref{corgal}.

\begin{Lemma} \label{lemegalbis} 
Let $M $,   $e_1,\ldots,e_k $, $C_1,\ldots,C_p$ be as in Lemma \ref{lemegal}. Let $\pi : \R^k \to \R^t$ be a $\R$-linear map, with $t\geq 1$. Then the rank of $(\pi(C_1),\ldots,\pi(C_p))$ in $\R^t$ (seen as a $\Q$-vector space) is equal to the  dimension of the minimal subspace $F$ of $\R^p$, defined over $\Q$, which contains the image of $\psi \circ  ^t \!\!\pi$; here $\psi$ is the $\R$-linear map of the dual of $\R^k$ to $\R^p$ which maps the canonical basis to $(e_1,\ldots,e_k)$.
\end{Lemma} 

\Dem of Lemma \ref{lemegalbis}: Let $P$ be the matrix of $\pi$ with respect to canonical bases, and $M' = PM$. Applying Lemma \ref{lemegal} to $M'$ gives directly the result.

\bigskip

\Dem of the second assertion of Corollary \ref{corgal}: Let $F$ denote the  minimal subspace  of $\R^p$, defined over $\Q$, which contains  the image of $\psi \circ  ^t \!\!\pi$; then Lemma \ref{lemegalbis} yields $\dim F = s$.  Now $\rk( ^t  \pi) = \rk(\pi) = t$ and $\psi$ is injective because $e_1,\ldots,e_k$ are $\R$-linearly independent (using conclusion $(ii)$ of Theorem \ref{thgal}), so that $\Im(\psi \circ  ^t \!\!\pi)$ has dimension $t$. Since this subspace is contained in both $F$ and $\Span_\R(e_1,\ldots,e_k) = \Im \psi$, we have $\dim(F\cap  \Span_\R(e_1,\ldots,e_k))\geq t$. Now the first part of Corollary \ref{corgal} (deduced in \S \ref{secenonce} from Theorem \ref{thgal}) shows that Eq. \eqref{eqcorgal} holds when $t$ is replaced with this (possibly larger) dimension; therefore it holds with $t$. This concludes the proof of  the second assertion of Corollary \ref{corgal}.

\subsection{Reduction to the non-oscillatory case} \label{subsecosc}

In this subsection, we deduce the general case of Theorem \ref{thgal} from the special case where $\om_1 =\ldots=\om_k =0$; notice that in this case we have $\phi_j \not\equiv \frac{\pi}2 \bmod \pi$ for any $j \in \unk$, so that Eq. \eqref{eqhypasyth} reads $|L_n(e_j) | = Q_n^{-\tau_j+o(1)}$.  This special case will be proved in the following subsections, under the assumption that $Q_{n+1} = Q_n^{1+o(1)}$ (which is weaker than the assumption $Q_{n+1} = Q_n ^{1+O(1/n)}$   we make when $\om_1$, \ldots, $\om_k$ may be non-zero). 

\bigskip

Let $\om_1$, \ldots, $\om_k$,  $\pha_1$, \ldots, $\pha_k$, and $(Q_n)$ be as in Theorem \ref{thgal}, with  $Q_{n+1} = Q_n ^{1+O(1/n)}$. Since there are  infinitely many integers $n$  such that, for any $j \in \unk$, $n \om_j + \pha_j \not\equiv \frac{\pi}2 \bmod \pi$, Proposition 1 of \cite{SFoscillate} provides $\eps , \lambda > 0$ and an increasing function $\psi : \N \to \N$ such that $\lim_{n \to \infty} \frac{\psi(n)}{n} = \lambda $ and, for any $n$ and any $j \in \unk$, $|\cos(\psi(n)\om_j+\pha_j)| \geq \eps  $. Let $L'_n = L_{\psi(n)}$ and $Q'_n = Q_{\psi(n)}$ for any $n \geq 1$. 
Then we have $|L'_n(e_j)| = {Q'_n}^{-\tau_j+o(1)}$ because $|\cos(\psi(n)\om_j+\pha_j)| =  Q_{\psi(n)}^{o(1)}$. Let us check that $Q'_{n+1} = Q_n'^{1+o(1)}$; then the special case of  Theorem \ref{thgal} will apply to the sequences $(L'_n)_{n\geq 1}$ and  $(Q'_n)_{n\geq 1}$, with the same other parameters: this will conclude the proof.

Since $Q_{n+1} = Q_n ^{1+O(1/n)}$ there exists $M > 0$ such that, for any $n \geq 1$, $Q_{n+1} \leq Q_n^{1+M/n}$; this implies 
$$\log Q_{n+\ell} \leq  (1 + M/n)^\ell \log Q_n $$
for any $\ell \geq 0$. Letting $\delta_n = \psi(n+1)-\psi(n) \geq 1$, we have:
\begin{eqnarray*}
\log Q'_{n+1} = \log Q_{\psi(n)+\delta_n} \leq (1 + M/\psi(n))^{\delta_n}  \log  Q_{\psi(n)} &\leq&   \exp(M \delta_n / \psi(n)) \log Q_{\psi(n)} \\
&=& (1+o(1))\log Q_n'
\end{eqnarray*}
since $1+x\leq e^x$ and $\delta_n = o(n)$ since $ \psi(n)  = \lambda n  + o(n)$. This concludes the reduction to  the  case where $\om_1 =\ldots=\om_k =0$ and $Q_{n+1} = Q_n^{1+o(1)}$.

\subsection{Proof of $(ii)$}   \label{subsecdemii}

Let us come now to  the easiest part of Theorem \ref{thgal}, namely $(ii)$. We shall prove simultaneously that    $e_1,\ldots,e_k $ are  linearly independent in $\R^p$, and that $F \cap \Q^p = \{(0,\ldots,0)\}$ where $F = \Span_\R(e_1,\ldots,e_k)$. With this aim in mind, we assume (by contradiction) that there exist real numbers    $\lam_1, \ldots, \lam_k $,  not all zero, such that $\sum_{j=1}^k \lam_j e_j \in  \Q^p$; multiplying all $\lam_j$ by a common denominator of the coordinates, we may assume  $\sum_{j=1}^k \lam_j e_j \in  \Z^p$.  Then $\kappa_{n}  = L_{n}(\sum_{j=1}^k \lam_j e_j)  = \sum_{j=1}^k \lam_j L_{n}(e_j) $  is an integer for any $n \geq 1$. Now if $n$ is sufficiently large then $|\kappa_{n}| \leq \sum_{j=1}^k | \lam_j |\, | L_{n}(e_j) |  < 1$, so that $\kappa_{n}  = 0$. Let $j_0$ denote the largest integer $j$ such that $\lam_j \neq 0$. Then for any  $n$   sufficiently large,  the fact that $\kappa_{n}  = 0$ implies $|\lam_{j_0} L_{n}(e_{j_0})| = |  \sum_{j=1}^{j_0-1} \lam_j L_{n}(e_j) |$ so that 
$$|\lam_{j_0}| \leq  \sum_{j=1}^{j_0-1}| \lam_j | \frac{   | L_{n}(e_j) |  }{ | L_{n}(e_{j_0}) |  } \leq  \sum_{j=1}^{j_0-1}| \lam_j | Q_n^{\tau_{j_0}-\tau_j+o(1)}$$
as $n \to \infty$. Now the right handside tends to 0 as $n \to \infty$ because we have  assumed that $\tau_1 > \ldots > \tau_k$, so that $\lam_{j_0} = 0$: this contradicts the definition of $\lam_{j_0}$.

Therefore such real numbers  $\lam_1, \ldots, \lam_k $ cannot exist, and this concludes the proof of~$(ii)$.

\subsection{Proof that $(ii)$ and $(iii)$ imply $(i)$} \label{subsecdemi}

Before proceeding in \S \S \ref{subsecpropc}  and \ref{subsec33} to the proof of $(iii)$, which is the main part, we deduce $(i)$ from $(ii)$ and $(iii)$. Recall that the second statement of $(i)$ is equivalent to the first one (which we shall prove now) thanks to Lemma \ref{lemegal} proved in \S \ref{subsecrank}.

Let $F$ be a subspace of $\R^p$, defined over $\Q$, which contains $e_1$, \ldots, $e_k$. Letting $s=\dim F$, we have $s > k$ using $(ii)$. 
Assertion $(iii)$ yields, for any $\eps>0$  and  any $Q$   sufficiently large (in terms of $\eps$),  a subset $\nvcalC$ and a lattice $\nvLam$ such that $\nvcalC \cap \nvLam = \{(0,\ldots,0)\}$.   Now $\nvcalC \cap F$   is a convex body, compact and symmetric with respect to the origin, in the Euclidean space $F$. On the other hand,  $\nvLam \cap F$ is a lattice in $F$ because $F$ is defined over $\Q$. Therefore  Minkowski's convex body theorem  (see for instance Chapter III of \cite{Cassels}) implies that $\nvcalC \cap F $ has volume less  than   $2^s \det(  \nvLam \cap F)$. 
Letting 
$$\alpha = \tau_1 + \ldots + \tau_k - \sigma_1 - \ldots - \sigma_{s-k} - s \eps,$$
 this volume is greater than or equal to $Q^\alpha$,   up to a multiplicative constant which depends only on $F$, $e_1 , \ldots,  e_k$,  $v_1$, \ldots, $v_p$ (using the inequalities $\sigma_1\geq \ldots \geq \sigma_p$). On the other hand, since $d_1\leq\ldots \leq d_p$ 
 we have $\det (\nvLam\cap F) \leq c Q^{\beta+o(1)}$ where $\beta = -d_1-\ldots - d_s$ and $c$ is a constant depending only on $F$.  
Since $Q$  can be chosen arbitrarily large, the above-mentioned consequence of  Minkowski's   theorem  yields $\alpha\leq \beta$. 
  Now $\eps$ can be any positive real number, so that we obtain
  $$ \tau_1 +\ldots + \tau_k + d_1 + \ldots + d_s \leq \sigma_1 + \ldots + \sigma_{s-k},$$
   thereby concluding the proof of $(i)$.

\subsection{A matrix lemma} \label{subsecpropc} 

We state and prove in this section the main tool in the proof of Theorem \ref{thintro2}, namely Lemma~\ref{lemmatpp}.
This result might be of independent interest; its proof relies on    estimating the  determinant and cofactors.

\begin{Lemma} \label{lemmatpp}
Let $A$ be a $k\times k$ matrix with real positive entries $a_{i,j}$, $1\leq i,j\leq k$, such that
 \begin{equation} \label{hypepsc}
 a_{i',j}a_{i,j'} \leq \frac1{(k+1)!}  a_{i,j} a_{i',j'} \mbox{ for any } i,j,i',j' \mbox{ such that } i<i' \mbox{ and } j<j'.
\end{equation}
Then $A$ is an invertible matrix, and letting $A^{-1} = [b_{i,j}]_{1\leq i,j\leq k}$ we have
$$|b_{j,i}| \leq \Big(1+ \frac1k + \frac1{k^2}\Big)  a_{i,j}^{-1} \mbox{ for any } i,j\in\unk.$$
\end{Lemma}

Lemma \ref{lemmatpp}  is optimal up to the value of the constant $1+ \frac1k + \frac1{k^2}$: it would be false with a constant less than $1/k$ instead (this is immediately seen by computing a diagonal coefficient of $AA^{-1}$, which is equal to 1).   We did not try to improve on the constant $1+ \frac1k + \frac1{k^2}$, but anyway it could easily be made smaller by replacing $ \frac1{(k+1)!} $ in \eqref{hypepsc} with a smaller constant.

\bigskip

In the proof of  Lemma \ref{lemmatpp}  we shall use the following lemma.

\begin{Lemma} \label{lemb}
Under the assumptions of Lemma \ref{lemmatpp}, for any $ \sigma \in \sk$    we have 
\begin{equation} \label{eqb1}
\prod_{j=1}^k  a_{\sigma(j),j}   \leq \etasig   \prod_{j=1}^k a_{j,j} 
\end{equation}
where $\etasig  = \frac1{(k+1)!} $ for $ \sigma \neq\Id$, and $\etaId = 1$.
\end{Lemma}

\Dem of Lemma \ref{lemb}: 
For $ \sigma \neq \Id$ let  $\capasig$ denote the largest integer $j \in \unk$ such that $\sigma(j) \neq j$; put also $\capaId = 0$. We are going to prove Eq. \eqref{eqb1}  by induction on $\capasig$.
If $\capasig \leq 1$ then $\sigma = \Id$, so that Eq. \eqref{eqb1} holds trivially. Let $\sigma \in \sk$ be such that $\capasig \geq 2$, and assume that  Eq. \eqref{eqb1} holds for any $\sigma'$ such that $\capasigpr < \capasig$. We have $\sigma(j)=j$ for any $j \in \capasigpuk$, and $\sigma(\capasig) < \capasig$. Let $j_0 = \sigma^{-1} (\capasig)$; then $j_0 < \capasig$. Let $\sigma' = \sigma \circ \tau_{j_0, \capasig}$ where $ \tau_{j_0, \capasig}$  is the transposition that exchanges $j_0$ and $ \capasig $. Then $\sigma'(j)=j$ for any $j \in \capasigk$ so that $\capasigpr < \capasig$ and Eq. \eqref{eqb1} holds for $\sigma'$. Since  $\sigma' (j) = \sigma(j) $ for $j \not\in \{j_0, \capasig\}$,  $\sigma'(j_0) = \sigma(\capasig)$ and $\sigma'(\capasig) = \capasig$, this implies (using the fact that $\etasigpr \leq 1$)
$$a_{\sigma(\capasig), j_0} a_{\capasig, \capasig} 
 \prod_{1 \leq j \leq k \atop  j \not\in \{j_0, \capasig\}} a_{\sigma(j),j}  \leq  \prod_{j=1}^k a_{j,j} .$$
On the other hand, Eq. \eqref{hypepsc} implies  
$$a_{ \capasig , j_0} a_{\sigma(\capasig), \capasig}  \leq  \frac1{(k+1)!}  a_{\sigma(\capasig), j_0} a_{\capasig, \capasig} $$
because  $\sigma(\capasig) < \capasig$ and $j_0 < \capasig $. Multiplying out the previous two inequalities yields Eq.~\eqref{eqb1} for $\sigma$, since $\sigma(j_0) = \capasig$. This concludes the proof of Lemma \ref{lemb}.

\bigskip

\Dem of Lemma \ref{lemmatpp}:  Letting $\Delta = | \det A\, |$ we have, using Lemma \ref{lemb}:
\begin{equation} \label{eqnvnva}
  \Delta  \geq   \prod_{j=1}^k a_{j,j}  - \sum_{\sigma \in \sk \atop \sigma \neq \Id} \prod_{j=1}^k  a_{\sigma(j),j}  \geq   \Big( 1 - \frac1{k+1}\Big)  \prod_{j=1}^k a_{j,j} > 0
 \end{equation}
 so that $A$ is invertible. Given $i,j\in\unk$ we have $|b_{j,i}| =\frac{\Delta_{i,j}}{\Delta}$ where $\Delta_{i,j}$ is the absolute value of the determinant of the matrix obtained from $A$ by deleting the $i$-th row and the $j$-th column. Using Lemma \ref{lemb} again we have
 \begin{equation} \label{eqnvnvb}
  \Delta_{i,j}  \leq \sum_{\sigma \in \sk \atop \sigma(j)=i }   \prod_{1 \leq j' \leq k \atop j' \neq j}  a_{\sigma(j'),j'} \leq \Big(  \sum_{\sigma \in \sk \atop \sigma(j)=i }  \etasig \Big) a_{i,j}^{-1}    \prod_{  j'=1}^k    a_{ j' ,j'}.
 \end{equation}
Now we have $\etasig = 1$ for at most one $\sigma$, and $\etasig = \frac1{(k+1)!}$ for all other permutations $\sigma$  among the $(k-1)!$ such that $\sigma(j)=i $, so that
$$\sum_{\sigma \in \sk \atop \sigma(j)=i }  \etasig \leq 1 + \frac{(k-1)!}{(k+1)!} = \frac{k+1+\frac1k}{k+1}.$$
Combining this upper bound with Eqns. \eqref{eqnvnva} and  \eqref{eqnvnvb}  yields
$$ |b_{j,i}| =\frac{\Delta_{i,j}}{\Delta} \leq  \frac{k+1+\frac1k}{k }a_{i,j}^{-1} ,$$
thereby completing the proof of Lemma \ref{lemmatpp}.

\subsection{Proof   of $(iii)$} \label{subsec33}

We are now in position to prove the remaining part of Theorem \ref{thintro2}, namely $(iii)$. We   assume  $\tau_1 > \ldots > \tau_k > 0$ and $\om_1 = \ldots = \om_k = 0$ (see \S \ref{subsecosc}), so that $|L_n(e_j) | = Q_n^{-\tau_j+o(1)}$. 

\bigskip

Before giving details, let us make a few comments on our strategy. 

Recall that Nesterenko's linear independence criterion is much easier to prove if the linear forms $L_n$, $L_{n+1}$, \ldots, $L_{n+p-1}$ are linearly independent (see  \S 2.3 of \cite{SFZu}    or the references to Siegel's criterion in  \S \ref{subsecsiegel} below). Of course this is not always the case, but Lemma \ref{lemmatpp} enables us to make  a step in this direction. Actually   letting $F = \Span_{\R}(e_1,\ldots,e_k)$, we consider the restrictions $L_{n | F}$ of the linear forms to $F$; recall that $\dim F = k$ thanks to $(ii)$ proved in \S \ref{subsecdemii}.  It is not true in general that $L_{n | F}$, $L_{n +1 | F}$, \ldots, $L_{n +k-1| F}$ are linearly independent linear forms on $F$: for instance, the equality $L_n = L_{n+1}$ might hold for any even integer $n$ (because of  the error terms $o(1)$ in the assumptions of  Theorem \ref{thgal}). To make this statement correct, we introduce a function $\ffi : \Netoile \to \Netoile$   such that $ \ffi(n) \geq n+1$ for any $n \geq 1$. The integer $\ffi(n)$ plays the role of $n+1$, that is: applying $\ffi$ corresponds to ``taking the next integer''. The idea is that $\ffi(n)$ will be large enough (in comparison to $n$) to avoid obvious counter-examples as above coming from error terms. In more precise terms, $ \ffi(n)$ will be defined by the property $Q_{\ffi(n)-1} \leq Q_n^{1+\nve} < Q_{\ffi(n)}$ (where $\nve$ is a  small positive real number); in this way, the error terms $o(1)$ in the assumptions of  Theorem \ref{thgal} will not be a problem any more.

With this definition, we shall prove that for any $n$ sufficiently large, the linear forms $L_{n | F}$, $L_{\ffi(n)  | F}$, $L_{\nvffi_2(n)  | F}$, \ldots, $L_{\nvffi_{k-1}(n)  | F}$ on $F$ are  linearly independent (where $\nvffi_i = \ffi \circ \ldots \circ \ffi$), so that they make up a basis of the dual vector space $F^\star$. In the proof of  Theorem \ref{thgal} we shall need the following quantitative version of this property: in writing the linear form $e_j^\star$ (defined by $e_j^\star (\lam_1 e_1+\ldots+\lam_k e_k)= \lam_j$) as a linear combination of $\frac1{L_n(e_j)} L_{n | F}$,    $\frac1{L_{\ffi(n)} (e_j)} L_{\ffi(n)  | F}$, \ldots,   $\frac1{L_{\nvffi_{k-1}(n)} (e_j)} L_{\nvffi_{k-1}(n)  | F}$, the coefficients that appear are bounded independently from $n$  (actually they are between $-3$ and $3$): see Eq. \eqref{eqcsqpp} below.   
This will follow from Lemma \ref{lemmatpp}  applied to the matrix $A_n = [|L_{\nvffi_{i-1}(n)} (e_j)|]_{1\leq i,j\leq k}$. The point in applying this lemma is that sharp upper and lower bounds on $|L_{\nvffi_{i-1}(n)} (e_j)|$ are available; the assumption $\tau_1> \ldots  > \tau_k$ plays also a central role here.

\bigskip

Now let us prove $(iii)$.

Let $\eps > 0$. We choose $\nve > 0$ sufficiently small, so that
\begin{equation} \label{eqdefnve}
((1+\nve)^{k-1}-1) \max(1, \tau_1 , \sigma_1) < \eps/4 .
\end{equation}
If $k=1$ there is no assumption on $\nve$, because it does not really appear in the proof:  Lemma \ref{lemmatpp} is a triviality in this case, and the proof of $(iii)$   reduces essentially   to that of \cite{SFZu}. 

\bigskip
 
For any $n \geq 1$, we define $\ffi(n)$ by $Q_{\ffi(n)-1} \leq Q_n^{1+\nve} < Q_{\ffi(n)}$, because the sequence $(Q_n)$ is increasing and we may assume $Q_n \geq 1$ for any $n$. Then we have $\ffi(n) \geq n+1$. This implies $\lim_{n \to + \infty} \ffi(n) = +\infty$, so that $Q_{\ffi(n)} = Q_{\ffi(n)-1} ^{1+o(1)}$  (because we assume $Q_{n+1} = Q_n ^{1+o(1)}$) and 
\begin{equation} \label{eqcarphi}
Q_{\ffi(n)} = Q_{n} ^{1+\nve + o(1)};
\end{equation}
here $o(1)$ denotes any sequence that tends to 0 as $n\to \infty$. Moreover   the assumption $|L_n(e_j)|   = Q_n^{-\tau_j+o(1)}$ implies  $|L_n(e_j) | > 0$  for any $ j$, $n$  with $n$ sufficiently large.  
We have also for any $n$ sufficiently large and any $j\in\unk$:
\begin{equation} \label{eq38bis}
|L_{\ffi(n)}(e_j)|  = Q_{\ffi(n)} ^{-\tau_j + o(1)} = Q_n^{-\tau_j(1 + \nve)+o(1)} < |L_n(e_j)|.
\end{equation}

For  $i \in \{0,\ldots, k-1\}$ let $\nvffi_i = \ffi \circ \ldots \circ \ffi$ denote the map $\ffi$ composed $i$ times with itself (so that $\nvffi_0(n) = n $ and $\nvffi_1(n) = \ffi(n)$). Put
$$A_n =  \Big[ |L_{\nvffi_{i-1}(n)} (e_j) |\Big]_{1\leq i,j\leq k}  $$
and denote by $a_{i,j}$ the entries of $A_n$ (omitting for simplicity the dependence on $n$).
Let us check the assumption \eqref{hypepsc} of  Lemma \ref{lemmatpp}, provided $n$ is  sufficiently large. Let $i,j,i',j'\in\unk$ be such that $i < i'$ and $j<j'$ ; we put $n' = \nvffi_{i-1}(n)$ and $n'' = \nvffi_{i'-1}(n)$, so that $n'' \geq \ffi(n')$. Using Eq. \eqref{eqcarphi} and the assumption $\tau_j > \tau_{j'}$ we obtain
$$\frac{a_{i',j} a_{i,j'}}{a_{i,j} a_{i',j'}} =
\Big| \frac{L_{n''}(e_j)  L_{n'}(e_{j'})  }{L_{n'}(e_j)   L_{n''}(e_{j'})  } \Big| 
= \frac{Q_{n''}^{\tau_{j'}-\tau_j+o(1)}}{Q_{n'}^{\tau_{j'}-\tau_j+o(1)}} \leq \Big(\frac{Q_{\ffi(n')}}{Q_{n'}}\Big)^{\tau_{j'}-\tau_j+o(1)} = Q_{n'}^{\nve(\tau_{j'}-\tau_j)+o(1)} \leq \frac1{(k+1)!}$$
if $n$ is sufficiently large, so that Lemma \ref{lemmatpp} applies.
Given $M = \sum_{j=1}^k \lambda_j e_j$ with $\lambda_1,\ldots,\lambda_k\in\R$, we have $L_{\nvffi_{i-1}(n)} (M) = \sum_{j=1}^k a_{i,j}\lambda'_j$ where we let $\lambda'_j = \lambda_j$ if $L_{\nvffi_{i-1}(n)} (e_j) > 0$, and $\lambda'_j = -\lambda_j$  otherwise. Therefore Lemma \ref{lemmatpp} yields, for any $j\in\unk$ and any $n$ sufficiently large:
\begin{equation}\label{eqcsqpp}
| \lambda_j | = |\lambda'_j | = \Big| \sum_{i=1}^k b_{j,i} L_{\nvffi_{i-1}(n)} (M) \Big| \leq \Big( 1 + \frac1k + \frac1{k^2}\Big) \sum_{i=1}^k \frac{ | L_{\nvffi_{i-1}(n)} (M)|}{| L_{\nvffi_{i-1}(n)} (e_j)|}.
\end{equation}
This upper bound on $|\lambda_j|$ in terms of the $| L_{\nvffi_{i-1}(n)} (M)| $ is the main tool we shall use now in the proof.

\bigskip

Let  $Q$ be sufficiently large in terms of $\eps$, and assume that $\nvcalC\cap\nvLam$ contains a non-zero point $P$. Then we have
$$P = \lam_1 e_1 + \ldots +  \lam_k e_k + u = (x_1,\ldots,x_p) \neq (0,\ldots,0)$$
with  $\lam_1,\ldots, \lam_k \in \R$, $u = \mu_1 v_1 + \ldots + \mu_p v_p  \in (\Span_\R(e_1,\ldots,e_k))^\perp$,  $| \lam_j | \leq Q^{\tau_j-\eps}$ for any $j \in\unk$, $|\mu_i| \leq Q^{-\sigma_i-\eps}$ for any $i\in\unp$, and  $\delta_{i,n} x_i \in \Z$ for any $i$, where   $n = \Psi(Q)$ is the largest integer such that $Q_n \leq Q$. In particular we have $Q_n \leq Q < Q_{n+1}$ so that $Q = Q_n^{1+o(1)}$, and $n$ tends to $\infty$ as $Q\to\infty$:  if $u_n = o(1)$, that is $u_n \to 0$ as $n\to\infty$, then $u_n$ tends also to 0 as $Q \to \infty$.

\bigskip

Let $\ell$ denote the least integer   such that
\begin{equation} \label{eqell}
\mbox{ for any   } j \in \unk,  \mbox{ we have }  |\lam_j L_{\ell}(e_j)| \leq \frac{\delta_{p,\ell}}{3k \delta_{p,n} } .
\end{equation}
Since $|\lambda_j| \leq Q^{\tau_j-\eps}$ and $n$ is sufficiently large, this upper bound holds for    $    n$ so that  this integer exists and we have $\ell \leq n$. 

\bigskip

The integer $\ell$ depends on $Q$ and on the choice of  a non-zero point $P \in \nvcalC\cap\nvLam$.  Let us prove that $\ell \to \infty$ as $Q \to\infty$, uniformly with respect to  the choice of $P$. Let $\ell_0\geq 1$, and denote by  $K_{\ell_0}  $   the set of all points $P'=  \lam'_1 e_1 + \ldots +  \lam'_k e_k+u'$ with 
$$|\lam'_j|  \min_{1 \leq \ell' \leq \ell_0} |L_{\ell'}(e_j)|  \leq \frac1{3k}  \mbox{ for any } j \in \unk,$$
 where $u' \in (\Span (e_1,\ldots,e_k))^\perp$ can be written as $u' = \mu'_1 v_1 + \ldots + \mu'_p v_p$ with $| \mu'_i | \leq Q^{-\sigma_i-\eps}$ for any $i\in\unp$.  By definition of $\ell$ and $K_{\ell_0}$, if $\ell \leq \ell_0$ then   $\frac{\delta_{p,n}}{\delta_{p,\ell}}P \in K_{\ell_0}$. Moreover the point  $\frac{\delta_{p,n}}{\delta_{p,\ell}}P$ belongs also to $\nvLamlz $ since
 $$\delta_{i,\ell_0} \Big( \frac{\delta_{p,n}}{\delta_{p,\ell}} x_i\Big) =  \Big(\frac{\delta_{i,\ell_0}}{\delta_{i,\ell}}\Big)\Big(\frac{\delta_{p,n}/\delta_{i,n}}{\delta_{p,\ell}/\delta_{i,\ell}} \Big)\Big( \delta_{i,n} x_i\Big) \in\Z$$
for any $i\in\unp$, by assumption on the divisors $\delta_{t,n}$. Therefore (assuming $\ell \leq \ell_0$) the point  $\frac{\delta_{p,n}}{\delta_{p,\ell}}P$ belongs to $ K_{\ell_0} \cap \nvLamlz $, which is a finite set because $K_{\ell_0}$ is compact and $\nvLamlz $ is discrete. Now   the function $\chi : K_{\ell_0} \cap \nvLamlz \to\R $ defined by $\chi(P')  = \norm{\pi_\perp(P')}$, where $\pi_\perp$ is the orthogonal projection on $ (\Span (e_1,\ldots,e_k))^\perp$, has a least positive value $\chi_0$. We have $\chi(\frac{\delta_{p,n}}{\delta_{p,\ell}}P) \neq 0$ because $P\not\in \Q^p\cap \Span (e_1,\ldots,e_k) = \{(0,\ldots,0)\}$ (using assertion $(ii)$ proved in \S \ref{subsecdemii}), so that 
$$\chi_0 \leq \chi \Big( \frac{\delta_{p,n}}{\delta_{p,\ell}}P \Big) =  \frac{\delta_{p,n}}{\delta_{p,\ell}} \norm{u} \leq Q_n^{d_p+o(1)} Q^{-\sigma_p-\eps} = Q^{d_p-\sigma_p-\eps+o(1)}$$
since $\delta_{p,\ell} \geq 1$ and $\sigma_p\leq \ldots \leq \sigma_1$. This inequality implies that  $Q$ is not too large in terms of $\ell_0$ and $\eps$ (because we assume $d_p \leq \sigma_p$). This   concludes the proof that $\ell \to \infty$ as $Q \to\infty$. In what follows, a sequence denoted by $o(1)$ will tend to 0 as $n$, $\ell$ or $ Q$ tends to $\infty$; therefore in any case, it tends to 0 as $ Q \to\infty$. Moreover, we may assume $\ell$ to be arbitrarily large.

\bigskip

We come back now to the point $P \in \nvcalC\cap\nvLam$ chosen above.
Since  $u = \mu_1 v_1 + \ldots + \mu_p v_p$  with $|\mu_h| \leq Q^{-\sigma_h-\eps}$ for any $h$, we have for any $i\in\unk$:
\begin{eqnarray}
|  L_{\nvffi_{i-1}(\ell)}(u)| &\leq&   \sum_{h=1}^p |\mu_h | |L_{\nvffi_{i-1}(\ell)}(v_h) |  \leq \sum_{h=1}^p Q^{-\sigma_h-\eps}  Q_{ \nvffi_{i-1}(\ell)}^{\sigma_h+o(1)} \nonumber \\
&\leq&  \sum_{h=1}^p  Q_n^{-\sigma_h-\eps+o(1)}  Q_{\ell}^{\sigma_h(1+\nve)^{i-1} +o(1 )} \mbox{ using Eq. \eqref{eqcarphi}}\nonumber\\
&\leq&    \sum_{h=1}^p \Big(\frac{Q_\ell}{ Q_n} \Big) ^{\sigma_h}  Q_n^{ -\eps+o(1)}   Q_{\ell}^{\eps/4+o(1 )}  \mbox{ using Eq. \eqref{eqdefnve} and the inequality $\sigma_h \leq \sigma_1$}\nonumber\\
&\leq&    \Big(\frac{Q_\ell}{ Q_n} \Big) ^{d_p} Q^{-\eps/2} < \frac13  \frac{\delta_{p,\ell}}{  \delta_{p,n} }    \mbox{ since $\sigma_h \geq \sigma_p \geq d_p $ and $\ell \leq n$}. \label{eq39bis}
\end{eqnarray} 

On the other hand, Eqns. \eqref{eqell} and \eqref{eq38bis} yield for any $i \in \unk$: 
$$|  L_{\nvffi_{i-1}(\ell)}(\sum_{j=1}^k \lam_j e_j)| \leq \sum_{j=1}^k  \frac{ |  L_{\nvffi_{i-1}(\ell)}(e_j)|}{ |  L_{ \ell }(e_j)|} \frac{\delta_{p,\ell}}{3k\delta_{p,n}} \leq \frac{\delta_{p,\ell}}{3 \delta_{p,n}},$$
since $\ell$ is sufficiently large. Combining this inequality with Eq. \eqref{eq39bis}
  we obtain for the point $P = \lam_1 e_1 + \ldots +  \lam_k e_k + u  $:
\begin{equation} \label{eqtroppetit}
  |  L_{\nvffi_{i-1}(\ell)}(P) | \leq   \frac{\delta_{p,\ell}}{3 \delta_{p,n}} + \frac{\delta_{p,\ell}}{3 \delta_{p,n}} < \frac{\delta_{p,\ell}}{  \delta_{p,n}} .
  \end{equation} 
Now we have $  L_{\nvffi_{i-1}(\ell)} = \ell_{1,\nvffi_{i-1}(\ell)}  X_1 + \ldots +  \ell_{p,\nvffi_{i-1}(\ell)}  X_p$ where $\ell_{j,\nvffi_{i-1}(\ell)} $ is a multiple of $\delta_{j,\nvffi_{i-1}(\ell)} $, and therefore of $\delta_{j,\ell}$ since $\nvffi_{i-1}(\ell)  \geq \ell$. 
Moreover  $\delta_{j,n} x_j\in\Z$ so that 
$$ \frac{\delta_{p,n}}{  \delta_{p,\ell} }  \ell_{j,\nvffi_{i-1}(\ell)} x_j = \Big( \frac{ \delta_{p,n} / \delta_{j,n}}{ \delta_{p,\ell} / \delta_{j,\ell}} \Big) \Big(\frac{\ell_{j,\nvffi_{i-1}(\ell)}}{\delta_{j,\ell}}\Big)(\delta_{j,n} x_j)\in\Z$$
since $\ell \leq n$, by assumption on the divisors $\delta_{t,n}$. Therefore we have $L_{\nvffi_{i-1}(\ell)}(P) \in  \frac{\delta_{p,\ell}}{ \delta_{p,n}}\Z$, and the upper bound \eqref{eqtroppetit} implies that this rational number is zero for any $i\in\unk$. Using  Eq. \eqref{eq39bis}
this yields the following  upper bound on  $|L_{\nvffi_{i-1}(\ell)}(M)|$ (where we let $M =  \sum_{j=1}^k \lam_j e_j$):
$$
 | L_{\nvffi_{i-1}(\ell)}(M) | =  | L_{\nvffi_{i-1}(\ell)}(u ) | \leq  \Big( \frac{Q_\ell}{Q_n}\Big)^{d_p} Q^{ -\eps/2}.
$$
Combining this upper bound with Eq. \eqref{eqcsqpp} yields, for any $j \in \unk$:
\begin{eqnarray*}
|\lam_j L_{\ell-1}(e_j) |
&\leq&  \Big( 1+\frac1k + \frac1{k^2} \Big)  \sum_{i=1}^k    \Big( \frac{Q_\ell}{Q_n}\Big)^{d_p} Q^{-\eps/2} Q_{ \nvffi_{i-1}(\ell)}^{ \tau_j+o(1)}Q_{\ell-1}^{-\tau_j+o(1)} \\
&\leq&  Q_\ell^{d_p+ \tau_j((1+\nve)^{i-1}-1)+o(1)} Q_n^{-d_p} Q^{-\eps/2}  \mbox{ using Eq. \eqref{eqcarphi}}\\
&\leq&  Q_\ell^{d_p+\eps / 4 + o(1)} Q_n^{-d_p} Q^{-\eps/2}  \mbox{ using the assumption $\tau_j \leq \tau_1$ and Eq. \eqref{eqdefnve}}\\
&\leq&  \Big( \frac{Q_\ell}{Q_n}\Big)^{d_p}  Q^{-\eps / 4 + o(1)} \leq  \frac{\delta_{p,\ell}}{ 3k \delta_{p,n}} \mbox{  since } Q_\ell \leq Q_n  = Q^{1+o(1)} \mbox{ and }  \frac{\delta_{p,\ell}}{  \delta_{p,n}} = \frac{Q_\ell^{d_p+o(1)}}{Q_n^{d_p+o(1)}}.
\end{eqnarray*}  
This contradicts the minimality of $\ell$ in Eq. \eqref{eqell}, thereby concluding the proof of $(iii)$.

\section{Consequences and related results} \label{sec2}

In this section we state and prove consequences of our main result (\S \S \ref{subsec21} and \ref{subsec22}), and mention Diophantine applications (\S \ref{secdio}). We also prove in \S \ref{subsecsiegel} an analogous result, in the spirit of Siegel's linear independence criterion. 

Throughout this section we restrict to the setting of Theorem \ref{thintro2}, omitting for simplicity the refinements of Theorem \ref{thgal} (eventhough they could have been adapted here).

\subsection{Distance to integers} \label{subsec21}

In this section we state corollaries of our criterion dealing with linear forms which are close to integers (rather than close to 0), as in Khintchine-Groshev's theorem for instance. In particular we deduce from Theorem \ref{thgal} a result (namely Corollary \ref{cortype2}  below) analogous to Nesterenko's linear independence criterion but which applies to sequences of simultaneous approximations of real numbers with the same denominator. This result is related to type  II  Pad\'e approximation  problems, in the same way as Nesterenko's  criterion  is related to  type~I     problems. In this respect, Theorem \ref{thgal} makes a bridge between the latter   and the former: it is related to Pad\'e approximation problems intermediate between  type  I   and  type~II  (see for instance \cite{Sorokinpi}). 

\bigskip

To begin with, let us state Theorem \ref{thintro2} in a {\em dual} way, namely in  terms of $C_1,\ldots,C_p\in \R^k$ rather than $e_1,\ldots,e_k\in\R^p$.

\begin{Th}\label{thcolumns}
Let $C_1,\ldots,C_p\in \R^k$, with $k,p\geq 1$. 

Let $\tau_1,\ldots,\tau_k $  and $(Q_n)_{n \geq 1}$  be as in Theorem \ref{thintro2}.

For any $n\geq 1$, let $\ell_{1,n},\ldots,\ell_{p,n}\in \Z$ be such that, as $n\to \infty$:
\begin{equation} \label{eqhypthintro1}
\max_{1 \leq i \leq p} |\ell_{i,n} | \leq Q_n^{1+o(1)} \mbox{\hspace{0.3cm} and  \hspace{0.3cm} }
\ell_{1,n} C_1 + \ldots + \ell_{p,n}C_p = 
\left( \begin{array}{c}
\pm Q_n^{-\tau_1+o(1)}   \\
\vdots\\
\pm Q_n^{-\tau_k+o(1)}  
\end{array} \right)
\end{equation}
where   the $\pm$ signs can be independent from one another. Then:
\begin{itemize}
\item[$(i)$] The rank of the family of vectors $C_1,\ldots,C_p$ in $\R^k$, considered as a $\Q$-vector space, is greater than or equal to 
$  k + \tau_1 +\ldots + \tau_k .$
\item[$(ii)$] For any non-zero linear form $\chi : \R^k \to \R$ there exists $i\in \unp$ such that $\chi(C_i) \not\in\Q$.
\item[$(iii)$] Let $\eps>0$, and $Q$ be sufficiently large  in terms of $\eps$. Let  $\lambda_1, \ldots, \lam_k \in \R$, not all zero, be such that $|\lam_j| \leq Q^{\tau_j-\eps}$ for any $j\in \unk$. Then denoting by $\chi$ the linear map $\R^k \to \R$ defined by $\chi(x_1,\ldots,x_k) = \lam_1x_1 + \ldots+\lam_kx_k$, we have
$${\rm dist}\Big(( \chi(C_1),\ldots,\chi(C_p)), \Z^p\setminus\{(0,\ldots,0)\}\Big)  \geq Q^{-1-\eps}$$
where  ${\rm dist}( y , \Z^p\setminus\{(0,\ldots,0)\})$ is the minimal distance of $y \in \R^p$ to a non-zero integer point.
\end{itemize}
\end{Th}

This result is just a translation of Theorem \ref{thintro2}. Indeed let us consider the matrix $M \in  {\rm Mat}_{k,p}(\R)$ of which $C_1,\ldots,C_p$ are the columns. We denote by $e_1,\ldots,e_k \in \R^p$ the rows of $M$. Then assumption \eqref{eqhypthintro1} means that the linear form $L_n = \ell_{1,n}X_1 + \ldots + \ell_{p,n}X_p$ on $\R^p$ is small at the points  $e_1,\ldots,e_k $. It is not difficult to see that $(ii)$ and $(iii)$ in  Theorem \ref{thcolumns} are respectively equivalent to  $(ii)$ and $(iii)$ in  Theorem  \ref{thintro2}, because $(\chi(C_1),\ldots,\chi(C_p)) = \lam_1e_1+\ldots+\lam_ke_k$. We remark also that assuming $k\leq p-1$ in Theorem \ref{thintro2} is not necessary; it has not been  used in the proof. This upper bound follows from $(ii)$, so that it is actually a consequence of the other assumptions.

\bigskip

Let us focus now on an important special case of Theorem \ref{thcolumns}, related to Pad\'e approximation: when $C_1$, \ldots, $C_k$ is the canonical basis of $\R^k$. This happens in all practical situations mentioned in \S \ref{secdio} below: indeed Pad\'e approximation provides linear combinations of $C_{k+1},\ldots,C_p$ which are very close to $\Z^k$.
 In this case, in $(ii)$ the interesting point is when the linear form  $\chi(x_1,\ldots,x_k) = \lam_1x_1 + \ldots+\lam_kx_k$ has rational coefficients $\lam_j$; then we have $\chi(C_i)\not\in \Q$ for some $i \in \{k+1,\ldots,p\}$. An analogous remark holds for $(iii)$; both are more easily stated as follows,  in terms of $e_1$,\ldots, $e_k$. We denote by $\norm{\cdot }$ any fixed norm on $\R^{p-k}$.

\begin{Cor} \label{coreprime} Under the assumptions  of Theorem \ref{thintro2}, suppose that for any $j\in \unk$ we have $e_j=(0,\ldots,0,1,0,\ldots,0,e'_j)$ with $e'_j\in \R^{p-k}$, where the 1 is in $j$-th position.

Then no non-trivial $\Q$-linear combination of $e'_1,\ldots,e'_k$ belongs to $\Q^{p-k}$. In addition, let $\eps>0$, and $Q$ be sufficiently large  in terms of $\eps$. Let  $\lambda_1, \ldots, \lam_k \in \Z$, not all zero, be such that $|\lam_j| \leq Q^{\tau_j-\eps}$ for any $j\in \unk$. Then for any $S \in \Z^{p-k}$ we have
 $$\norm{\lam_1e'_1+\ldots+\lam_ke'_k-S}\geq Q^{-1-\eps}.$$
\end{Cor}

This corollary is a measure of linear independence of the vectors $e'_1,\ldots,e'_k$ and those of the canonical basis of $\Z^{p-k}$. It can be weakened by assuming $|\lam_j| \leq Q^{\tau-\eps}$ for any $j\in \unk$, where $\tau = \min(\tau_1,\ldots,\tau_k)$ (as in Theorem \ref{thdist} below). Then a {\em measure of non-discreteness} (in the sense of \cite{Gutnik2003}) is obtained for the lattice $\Z e'_1+ \ldots + \Z e'_k + \Z^{p-k}$, which has rank $p$. In the examples  \eqref{eqGu1}, \eqref{eqGu2} and  \eqref{eqHP} considered in \S \ref{secdio} below, the matrix with columns $C_{k+1}$, \ldots, $C_p$ is symmetric (with $p = 2k$), so that this lattice is exactly $\Z C_1 + \ldots + \Z C_p$ (using  the fact that $C_1$, \ldots, $C_k$ is the canonical basis of $\R^k$).

This case $k = p/2$ lies ``in the middle'' between $k=1$, which corresponds to type  I  Pad\'e approximation and Nesterenko's original criterion, and $k=p-1$, which corresponds to type  II Pad\'e approximation. In the latter case, Corollary \ref{coreprime} yields the following result by letting $\xi_j = -e'_j$.

\begin{Cor} \label{cortype2} 
Let $k \geq 1$, and $\xi_1,\ldots,\xi_k\in \R$.

Let $\tau_1,\ldots,\tau_k >0 $ be pairwise distinct  real numbers. 

Let $(Q_n)_{n \geq 1}$  be an increasing sequence of positive integers, such that   
$Q_{n+1} = Q_n^{1+o(1)}$.

For any $n\geq 1$, let $\ell_{1,n}, \ldots,  \ell_{k,n},  \ell_{k+1,n}  \in \Z$ be such that 
$$\max_{1 \leq i \leq k+1} |\ell_{i,n} | \leq Q_n^{1+o(1)} $$
and
$$|\ell_{k+1,n} \xi_j - \ell_{j,n} | = Q_n^{-\tau_j+o(1)}\mbox{ for any } j\in \unk.$$
Then:
\begin{enumerate}
\item[$(i)$] The numbers 1, $\xi_1,\ldots,\xi_k$ are $\Q$-linearly independent.
\item[$(ii)$]  Let $\eps>0$, and $Q$ be sufficiently large (in terms of $\eps$). 
Then for any $(a_0,a_1,\ldots,a_k)\in\Z^{k+1}\setminus\{(0,\ldots,0)\}$ with $|a_j| \leq Q^{\tau_j-\eps}$ for any $j\in\unk$, we have:
$$|a_0+a_1\xi_1+\ldots+a_k\xi_k|\geq Q^{-1-\eps}.$$
\end{enumerate}
\end{Cor}

We have not found this statement  in the literature;  see however \cite{EMS} (p. 98),  \cite{Hata1993AA} (Lemma 2.1) or \cite{Hata1998TransAMS} (Lemma 6.1) for  related results, which are probably closer to Siegel's criterion than to Nesterenko's (see \S \ref{subsecsiegel} below).

\subsection{Upper bound on a Diophantine exponent} \label{subsec22}

Given a subspace $F$ of $\R^p$, and a non-zero point $P \in \R^p$, we denote by $\Dist (P, F )$ the projective distance of $P$ to $F$, seen in ${\mathbb P}^p(\R)$. Several definitions may be given, all of them equivalent up to multiplicative constants (see for instance \cite{SchmidtAnnals}); we choose $\Dist (P, F ) = \frac{\norm{u}}{\norm{P}}$ where $u$ is the orthogonal projection of $P$ on $F^{\perp}$ (that is, $P$ can be written as $u+f$ with $ u \in F^{\perp}$  and $f\in F$), and $\norm{\cdot} $ is the Euclidean norm on $\R^p$.

\smallskip

The following result is a consequence of Theorem \ref{thintro2}.

\begin{Th} \label{thdist}
Under the assumptions of Theorem \ref{thintro2}, let $\tau = \min(\tau_1,\ldots, \tau_k)$ and $F = \Span_\R(e_1,\ldots,e_k)$. Then for any $\eps > 0$ and any $P \in \Z^p \moins \{(0,\ldots,0)\}$ we have:
$$\Dist (P, F ) \geq \norm{P}^{-1-\frac{1}{\tau}-\eps}$$
provided $\norm{P}$ is sufficiently large in terms of $\eps$.
\end{Th}

It is important to notice that Theorem \ref{thdist} is {\em not} optimal, since it involves only $\min(\tau_1,\ldots, \tau_k)$. It is specially interesting when $\tau_1,\ldots, \tau_k$ are close to one another.

The interest of Theorem \ref{thdist} is that it can be written as an upper bound on a Diophantine exponent which measures the approximation of $F$ by points of $\Z^p$ (see \cite{SchmidtAnnals}, \cite{omegadun}, \cite{omegadde}).

\bigskip

\Dem of Theorem \ref{thdist}: Using assertion $(ii)$ of Theorem \ref{thintro2}, we see that $(e_1,\ldots, e_k)$ is a basis of $F$. Since $F$ is finite-dimensional, all norms on $F$ are equivalent: there exists $\kappa > 0$ such that, for any $f = \lam_1 e_1 + \ldots + \lam_k e_k \in F$ (with $\lam_j \in \R$), we have $\max|\lam_j| \leq \kappa \norm{f}$.

Let $\eps>0$ be such that $\eps < \tau$. Let $Q_0$ be such that assertion $(iii)$ of Theorem \ref{thintro2} holds for any $Q \geq Q_0$; we assume that $\norm{P} \geq Q_0^{\tau-\eps} / \kappa$. Letting $Q  = (\kappa \norm{P})^{1/(\tau-\eps)}$ we have $Q \geq Q_0$. Since $P \in  \Z^p \moins \{(0,\ldots,0)\}$, $P$ does not belong to the set $\nvcalC$ defined in assertion $(iii)$. Now writing $P = \lam_1 e_1 + \ldots + \lam_k e_k +u$ with $\lam_j \in \R$ and $u \in  F^\perp$, we have 
$$\max_{1\leq j \leq k} |\lam_j| \leq \kappa \norm{  \lam_1 e_1 + \ldots + \lam_k e_k } \leq  \kappa \norm{P} = Q^{\tau-\eps}$$
so that $\norm{u} > Q^{-1-\eps}$. Using the definition of $Q$ and that of $\Dist (P, F )$, this concludes the proof of Theorem \ref{thdist}.

\subsection{Connection with a Siegel-type criterion} \label{subsecsiegel}

The following result is analogous to Theorem \ref{thintro2}, but its proof is much easier. It relies on Siegel's ideas for linear independence (see for instance \cite{EMS}, p. 81--82 and 215--216,  or  \cite{Marcovecchio}, Proposition 4.1).  Special cases of this result have already been used in Diophantine results (see \S \ref{secdio} below).

\begin{Prop}\label{propsiegel}
Let $1\leq k \leq p-1$,  and   $e_1,\ldots,e_k \in \R^p$ be $\R$-linearly independent vectors. 

Let $(Q_n)_{n \geq 1}$  be an increasing sequence of positive integers, and for 
 any $n\geq 1$, let $ L_n^{(t)} = \ell_{1,n}^{(t)}X_1 + \ldots + \ell_{p,n}^{(t)}X_p$  be $p$ linearly independent  linear forms on $\R^p$ (for $1\leq t \leq p$), with integer coefficients $\ell_{i,n}^{(t)}$ such  that, as $n\to\infty$:
$$ |L_n^{(t)} (e_j) | \leq Q_n^{-\tau_j+o(1)}    \mbox{ for any } j \in \unk    \mbox{ and any } t\in \unp,
$$
where $\tau_1,\ldots,\tau_k >0 $ are  real numbers, and 
$$\max_{1 \leq i \leq p \atop 1\leq t \leq p} |\ell_{i,n}^{(t)}  | \leq Q_n^{1+o(1)}.$$
Then:
\begin{enumerate}
\item[$(a)$] Conclusions $(i)$ and $(ii)$ of Theorem \ref{thintro2} hold.
\item[$(b)$]    Let $\eps>0$, and $n$ be sufficiently large (in terms of $\eps$). Let $\calC_n$ denote the set of all vectors that can be written as $\lambda_1 e_1 + \ldots + \lambda_k e_k  + u$ with:
$$\left\{ \begin{array}{l}
 \lambda_1,\ldots,\lambda_k\in \R  \mbox{  such that } |\lambda_j| \leq Q_n^{\tau_j-\eps}   \mbox{ for any } j\in\unk\\
 u \in (\Span_\R(e_1,\ldots,e_k))^\perp \mbox{ such that } \norm{u} \leq Q_n^{-1-\eps} 
\end{array}\right.
$$
Then $\calC_n \cap \Z^p = \{(0,\ldots,0)\}$.
\end{enumerate}
\end{Prop}

The main difference with  Theorem \ref{thintro2} is that we require here $p$  linearly independent  linear forms for any $n$ (and we also assume $e_1,\ldots,e_k $ to  be $\R$-linearly independent). This  makes the proof much easier, and enables one to get rid of several important assumptions of  Theorem \ref{thintro2} (namely $Q_{n+1} = Q_n^{1+o(1)}$, $\tau_1,\ldots,\tau_k   $ pairwise distinct, and $|L_n(e_j) | $ not too small). 

If $Q_{n+1} = Q_n^{1+o(1)}$ in Proposition \ref{propsiegel} then in $(b)$ we may replace $Q_n$ with any $Q$, by letting $n$ be such that $Q_n \leq Q < Q_{n+1}$.

\bigskip

\Dem of Proposition \ref{propsiegel}: To prove conclusion $(i)$ of Theorem \ref{thintro2}, let $F$ be a subspace of $\R^p$  defined over $\Q$, of dimension $d$,  which contains $e_1$, \ldots, $e_k$. Let $n$ be sufficiently large. Up to reordering 
$L_n^{(1)}$, \ldots, $L_n^{(p)}$, we may assume the restrictions of $L_n^{(1)}$, \ldots, $L_n^{(d)}$ to $F$ to be   linearly independent  linear forms on $F$. Denoting by $(u_1,\ldots,u_d)$ a basis of $F$ consisting in vectors of $\Z^p$, the matrix $[L_n^{(t)}(u_j)]_{1\leq t,j\leq d}$ has a non-zero integer determinant. By making suitable linear combinations of the columns, the values $L_n^{(t)}(e_1)$, \ldots,  $L_n^{(t)}(e_k)$ appear and lead to the upper bound $Q_n^{d-k-\tau_1 -\ldots - \tau_k+o(1)}$ on the absolute value of this determinant. This concludes the proof of  $(i)$ of Theorem \ref{thintro2}.

To prove part $(b)$ of Proposition \ref{propsiegel} (which implies conclusion $(ii)$ of Theorem \ref{thintro2}), we let $P = \lambda_1 e_1 + \ldots + \lambda_k e_k  + u \in \calC_n \cap \Z^p$ be non-zero; then $L_n^{(t)}(P)\neq 0$ for some $t$, but $L_n^{(t)}(P)\in\Z$ and $|L_n^{(t)}(P)|<1$. This concludes  the proof of Proposition \ref{propsiegel}.

\subsection{Diophantine applications} \label{secdio}

The main interest of Theorems \ref{thintro2} and \ref{thgal} is that they provide (in conclusion $(i)$) a lower bound for the rank of $(C_1,\ldots,C_p)$. Such a lower bound (with $k$ essentially equal to $a^\eps$) implies Theorem \ref{thdistrib}, using a general lemma of linear algebra (see \cite{SFdistrib} for details). This kind of lower bounds (with $k\geq 2$) exists in the literature: for instance  Gutnik proved \cite{Gutnik83} that the vectors
\begin{equation} \label{eqGu1}
\left(\begin{array}{c} 1 \\ 0 \end{array} \right), \hspace{0.3cm}
\left(\begin{array}{c} 0 \\ 1 \end{array} \right),  \hspace{0.3cm}
\left(\begin{array}{c} -2 \log 2 \\   \zeta(2) \end{array} \right),  \hspace{0.3cm}
\left(\begin{array}{c} \zeta(2) \\ -3 \zeta(3) \end{array} \right)
\end{equation}
are $\Q$-linearly independent in $\R^2$ (so that, for any $r\in \Q^\star$, at least one number among $\zeta(2) - 2r \log 2 $ and $3\zeta(3) - r\zeta(2)$ is irrational). More recently he obtained also \cite{Gutnik2003} the $\Q$-linear independence of
\begin{equation} \label{eqGu2}
\left(\begin{array}{c} 1 \\ 0 \end{array} \right), \hspace{0.3cm}
\left(\begin{array}{c} 0 \\ 1 \end{array} \right),  \hspace{0.3cm}
\left(\begin{array}{c}  2 \zeta(3) \\   3\zeta(4) \end{array} \right),  \hspace{0.3cm}
\left(\begin{array}{c} 3 \zeta(4) \\  6 \zeta(5) \end{array} \right).
\end{equation}
In the same spirit,  T.~Hessami-Pilehrood has proved \cite{Hessami} that if $q$ is greater than some explicit function of $k$ then the following $2k$ vectors are $\Q$-linearly independent in $\R^k$:
\newcommand{\musq}{ {\tiny\frac{-1}{q }}}
\begin{equation} \label{eqHP}
\left(\begin{array}{c} 1 \\ 0 \\ \vdots \\ 0 \end{array} \right), \hspace{0.3cm}
\left(\begin{array}{c} 0 \\ 1 \\ \vdots \\ 0 \end{array} \right), \hspace{0.3cm}
\ldots, \hspace{0.3cm}
\left(\begin{array}{c} 0 \\ 0 \\ \vdots \\ 1 \end{array} \right), \hspace{0.3cm}
\end{equation}
$$
\left(\begin{array}{c} \Li_1(\musq)  \\  \Li_2(\musq) \\ \vdots \\  \Li_k(\musq) \end{array} \right), \hspace{0.3cm}
\ldots, \hspace{0.3cm}
\left(\begin{array}{c} \combitiny{j-1}{j-1} \Li_j(\musq)  \\  \combitiny{j}{j-1} \Li_{j+1}(\musq) \\ \vdots \\  \combitiny{j+k-2}{j-1} \Li_{j+k-1}(\musq) \end{array} \right), \hspace{0.3cm}
\ldots, \hspace{0.3cm}
\left(\begin{array}{c} \combitiny{k-1}{k-1} \Li_k(\musq)  \\  \combitiny{k}{k-1} \Li_{k+1}(\musq) \\ \vdots \\  \combitiny{2k-2}{k-1} \Li_{2k-1}(\musq) \end{array} \right)   .
$$
The same result holds with $1/q$ instead of $-1/q$; see also Gutnik's preprints cited in  \cite{Hessami}. 

\bigskip

These results share two common features: they rely on a special case of Proposition~\ref{propsiegel}, and they prove the linear independence of the full set of $p$ vectors involved. Using Theorem~\ref{thgal} it should not be difficult to produce alternative proofs of these results, in which only one sequence of small linear forms is constructed (instead of $p$ linearly independent ones); this may lead to further generalizations: for instance no proof of Ball-Rivoal's lower bound \eqref{eqBR} is known without using Nesterenko's criterion. Moreover, it should be possible also to obtain lower bounds for the rank of a family of vectors (like \eqref{eqGu1} or \eqref{eqGu2} up to $\zeta(a)$, or \eqref{eqHP} with smaller values of $q$)  eventhough the present methods fail to prove the linear independence of the full set.

\newcommand{\url}{\texttt}\providecommand{\bysame}{\leavevmode ---\ }
\providecommand{\og}{``}
\providecommand{\fg}{''}
\providecommand{\smfandname}{\&}
\providecommand{\smfedsname}{\'eds.}
\providecommand{\smfedname}{\'ed.}
\providecommand{\smfmastersthesisname}{M\'emoire}
\providecommand{\smfphdthesisname}{Th\`ese}


\begin{thebibliography}{10}

\bibitem{Apery}
{\scshape R.~Ap{\'e}ry} -- {\og Irrationalit\'e de $\zeta(2)$ et
  $\zeta(3)$\fg}, in \emph{Journ{\'e}es Arithm{\'e}tiques (Luminy, 1978)},
  Ast{\'e}risque, no.~61, 1979, p.~11--13.

\bibitem{BR}
{\scshape K.~Ball {\normalfont \smfandname} T.~Rivoal} -- {\og
  Irrationalit{\'e} d'une infinit{\'e} de valeurs de la fonction z\^eta aux
  entiers impairs\fg}, \emph{Invent. Math.} \textbf{146} (2001), no.~1,
  p.~193--207.

\bibitem{Boualg}
{\scshape N.~Bourbaki} -- {\og Alg\`ebre\fg}, ch.~II, Hermann, third.
  \smfedname, 1962.

\bibitem{omegadde}
{\scshape Y.~Bugeaud {\normalfont \smfandname} M.~Laurent} -- {\og On transfer
  inequalities in {D}iophantine approximation, {$II$}\fg}, \emph{Math. Z.}
  \textbf{265} (2010), p.~249--262.

\bibitem{Cassels}
{\scshape J.~Cassels} -- \emph{An introduction to the geometry of numbers},
  Grundlehren der Math. Wiss., no.~99, Springer, 1959.

\bibitem{EMS}
{\scshape N.~Fel'dman {\normalfont \smfandname} Y.~Nesterenko} -- \emph{Number
  theory {IV}, transcendental numbers}, Encyclopaedia of Mathematical Sciences,
  no.~44, Springer, 1998, A.N. Parshin and I.R. Shafarevich, eds.

\bibitem{SFoscillate}
{\scshape S.~Fischler} -- {\og Nesterenko's criterion when the small linear
  forms oscillate\fg}, \emph{Archiv der Math.} \textbf{98} (2012), no.~2,
  p.~143--151.

\bibitem{SFdistrib}
\bysame , {\og Distribution of irrational zeta values\fg}, manuscript, October 2013;  will   be
  posted soon on arxiv.

\bibitem{FHKL}
{\scshape S.~Fischler, M.~Hussain, S.~Kristensen {\normalfont \smfandname}
  J.~Levesley} -- {\og A converse to linear independence criteria, valid almost
  everywhere\fg}, preprint arxiv 1302.1952 [math.NT], submitted, 2013.

\bibitem{eddzero}
{\scshape S.~Fischler {\normalfont \smfandname} T.~Rivoal} -- {\og
  Irrationality exponent and rational approximations with prescribed
  growth\fg}, \emph{Proc. Amer. Math. Soc.} \textbf{138} (2010), no.~8,
  p.~799--808.

\bibitem{SFZu}
{\scshape S.~Fischler {\normalfont \smfandname} W.~Zudilin} -- {\og A
  refinement of {N}esterenko's linear independence criterion with applications
  to zeta values\fg}, \emph{Math. Ann.} \textbf{347} (2010), p.~739--763.

\bibitem{Gutnik83}
{\scshape L.~Gutnik} -- {\og On the irrationality of some quantities containing
  $\zeta(3)$\fg}, \emph{Acta Arith.} \textbf{42} (1983), no.~3, p.~255--264,
  (in Russian) ; translation in Amer. Math. Soc. Transl. {\bf 140} (1988), p.
  45--55.

\bibitem{Gutnik2003}
\bysame , {\og On linear forms with coefficients in {$\N \zeta(1 + \N)$}\fg},
  in \emph{Proceedings of the Session in analytic number theory and Diophantine
  equations (Bonn, 2002)} (D.~Heath-Brown {\normalfont \smfandname} B.~Moroz,
  \smfedsname), Bonner Mathematische Schriften, no. 360, 2003, p.~1--45.

\bibitem{Hata1993AA}
{\scshape M.~Hata} -- {\og Rational approximations to $\pi$ and some other
  numbers\fg}, \emph{Acta Arith.} \textbf{63} (1993), no.~4, p.~335--349.

\bibitem{Hata1998TransAMS}
\bysame , {\og The irrationality of $\log(1+1/q)\log(1-1/q)$\fg}, \emph{Trans.
  Amer. Math. Soc.} \textbf{350} (1998), no.~6, p.~2311--2327.

\bibitem{Hessami}
{\scshape T.~Hessami~Pilehrood} -- {\og Linear independence of vectors with
  polylogarithmic coordinates\fg}, \emph{Vestnik Moskov. Univ. Ser. I Mat.
  Mekh. [Moscow Univ. Math. Bull.]} \textbf{54} (1999), no.~6, p.~54--56
  [40--42].

\bibitem{omegadun}
{\scshape M.~Laurent} -- {\og On transfer inequalities in {D}iophantine
  approximation\fg}, in \emph{Analytic Number Theory, Essays in Honour of
  {K}laus {R}oth} (W.~Chen, W.~Gowers, H.~Halberstam, W.~Schmidt {\normalfont
  \smfandname} R.~Vaughan, \smfedsname), Cambridge Univ. Press, 2009,
  p.~306--314.

\bibitem{Marcovecchio}
{\scshape R.~Marcovecchio} -- {\og Linear independence of linear forms in
  polylogarithms\fg}, \emph{Annali Scuola Norm. Sup. Pisa} \textbf{V} (2006),
  no.~1, p.~1--11.

\bibitem{Nesterenkocritere}
{\scshape Y.~Nesterenko} -- {\og On the linear independence of numbers\fg},
  \emph{Vestnik Moskov. Univ. Ser. I Mat. Mekh. [Moscow Univ. Math. Bull.]}
  \textbf{40} (1985), no.~1, p.~46--49 [69--74].

\bibitem{RivoalCRAS}
{\scshape T.~Rivoal} -- {\og La fonction z\^eta de {R}iemann prend une
  infinit{\'e} de valeurs irrationnelles aux entiers impairs\fg}, \emph{C. R.
  Acad. Sci. Paris, Ser. I} \textbf{331} (2000), no.~4, p.~267--270.

\bibitem{SchmidtAnnals}
{\scshape W.~Schmidt} -- {\og On heights of algebraic subspaces and
  {D}iophantine approximations\fg}, \emph{Annals of Math.} \textbf{85} (1967),
  p.~430--472.

\bibitem{Sorokinpi}
{\scshape V.~Sorokin} -- {\og A transcendence measure for {$\pi^2$}\fg},
  \emph{Mat. Sbornik [Sb. Math.]} \textbf{187} (1996), no.~12, p.~87--120
  [1819--1852].

\end{thebibliography}
\end{document}